
\input amstex
\documentstyle{amsppt}
\nologo
\NoBlackBoxes
\mag=1200
\hsize=31 pc
\vsize=44 pc
\hcorrection{5mm}

\topmatter
\title   Immersed surfaces and Dehn surgery
\endtitle
\author  Ying-Qing Wu$^1$
\endauthor
\leftheadtext{Ying-Qing Wu}
\rightheadtext{Immersed surfaces and Dehn surgery}
\address Department of Mathematics, University of Iowa, Iowa City, IA
52242
\endaddress
\email  wu\@math.uiowa.edu
\endemail

\thanks  $^1$ Partially supported by NSF grant \#DMS 9802558
\endthanks


\endtopmatter

\document

\define\proof{\demo{Proof}}
\define\endproof{\qed \enddemo}
\define\a{\alpha}
\redefine\b{\beta}
\redefine\d{\delta}
\define\r{\gamma}

\define\area{\text{\rm Area}}

\redefine\e{\epsilon}
\redefine\bar{\overline}
\redefine\H{\Bbb H}
\redefine\Z{\Bbb Z}

\redefine\l{\lambda}
\redefine\bdd{\partial}
\define\Int{\text{\rm Int}}

\baselineskip 15pt
\input epsf.tex

\head \S 1.  Introduction
\endhead

The problem of how many Dehn filling on a torus boundary component $T$
of a 3-manifold $M$ will make a closed embedded essential surface $F$
compressible has been settled.  A slope $\b$ on $T$ is a coannular
slope if it is homotopic to some curve on $F$.  As an embedded
essential surface, $F$ can have at most one coannular slope.  If $F$
has a coannular slope $\b$ on $T$, then by a result of
Culler-Gordon-Luecke-Shalen [CGLS, Theorem 2.4.3], $F$ is
incompressible in all $M(\r)$ such that $\Delta(\b, \r) > 1$, where
$\Delta(\b, \r)$ denotes the minimal intersection number between the
slopes $\b$ and $\r$.  If $F$ has no coannular slopes, then it is
incompressible in the Dehn filling space $M(\r)$ for all but at most
three $\r$ [Wu].  There are examples showing that these are the best
possible.

While many manifolds do not contain embedded essential surfaces, it
has been shown by Cooper, Long and Reid [CLR] that most bounded
3-manifolds, in particular all hyperbolic manifolds with some toroidal
boundaries, contain immersed closed essential surfaces.  There has
been a lot of interest recently on immersed surfaces, see for example
[AR,CLR,CL1,CL2,Oe,Re].  It seems important to understand to which
extent the above theorems for embedded surfaces can be generalized to
immersed surfaces.

Let $S$ be a surface of finite type, i.e.\ compact surface with
finitely many points removed.  $S$ may be disconnected or
unorientable.  We define a {\it surface\/} (of type $S$) in $M$ to be
a continuous piecewise smooth map $F: S \to M$ which is an immersion
almost everywhere.  $F$ is {\it hyperbolic\/} if all components of $F$
have negative Euler characteristic.  A compact 3-manifold $M$ is
hyperbolic if its interior admits a complete hyperbolic structure.

Let $T$ be a set of tori in $\bdd M$.  A curve on a surface is {\it
simple\/} if it has no self-intersection.  A {\it slope\/} $\r_i$ on
$T$ is the isotopy class of a simple nontrivial curve on $T$.  A slope
$\r$ is a{\it coannular slope\/} if some nontrivial multiple of $\r$
is homotopic to a curve on $F$, in which case we say that $F$ is
coannular to $T$.  A {\it multiple slope\/} $\r = (\r_1, ..., \r_n)$
on $T$ is a set of slopes $\r_i$, one for each component $T_i$ of $T$.
Denote by $M(\r)$ the {\it Dehn filling space\/} along $\r$, i.e.\ the
manifold obtained by attaching a solid torus $V_i$ to each $T_i$
($i\leq n$) so that $\r_i$ bounds a meridian disk in $V_i$.  Given two
slopes $\a, \b$ on a torus, we use $\Delta(\a, \b)$ to denote the
minimal geometric intersection number between $\a$ and $\b$.  If $\r$
is a multiple slope and $\b$ is a finite set of slopes on $T$, denote
$$ \Delta(\r, \b) = \min\, \{ \, \Delta(\r_i, \b_j) \, | \, \b_j
\subset T_i \, \}$$ In particular, $\Delta(\r, \b) > 0$ if and only if
$\r_i \notin \b$ for all $i$.  Note that $\b$ may have none or
finitely many slopes on a component $T_i$ of $T$.  The following is
our main theorem.

\proclaim{Theorem 5.3} Let $T$ be a set of tori on the boundary of a
compact, orientable, hyperbolic 3-manifold $W$.  Let $F$ be a compact
essential surface in $W$ with $\bdd F \subset \bdd M - T$, and let
$\b$ be the set of coannular slopes of $F$ on $T$.  Then there is an
integer $K$ and a finite set of slopes $\Lambda$ on $T$, such that $F$
is $\pi_1$-injective in $W(\r)$ for all multiple slopes $\r$ on $T$
satisfying $\Delta(\r, \b) \geq K$ and $\r_i \notin \Lambda$.
\endproclaim

The result is best possible in the sense that there is no universal
bound on the constant $K$, see Theorem 6.1.  Note that $F$ has only
finitely many coannular slopes on $T$, i.e.\ $\b$ is a finite set.
See the remark before Theorem 5.2.  Thus in certain sense, Theorem 5.2
says that $F$ survives most Dehn fillings on $M$.  In particular, if
$F$ is not coannular to $T$, then $F$ survives all surgeries after
excluding a finite set of slopes on each component of $T$.

When $T$ has only one component, Theorem 5.3 can be generalized to
arbitrary compact orientable 3-manifolds $M$.  However the theorem is
no longer true when $M$ contains some Seifert fibred submanifolds and
$T$ contains more than one components.  An easy example is when $T$ is
a pair of tori $T_1, T_2$ coannular to each other.  If $F$ is
compressible in $M(\r_1, \r_2)$, then it is compressible in $M(\r_1',
\r'_2)$ for all $(\r'_1, \r'_2)$ obtained by twisting $(\r_1, \r_2)$
along an essential annulus with one boundary component on each of
$T_i$.  More complicated examples can be constructed where no two
components of $T$ are coannular.  However, the theorem is true if one
further excludes all slopes of distance at most one from the fiber
slopes.  More details will appear elsewhere.

Another interesting topic is to construct immersed essential surfaces
in hyperbolic 3-manifolds.  See for example [AR, CLR, FF, CL1, CL2,
Li].  One of the most important method is the Freedman tubing [FF].
Given a proper surface $F$ in $M$, a {\it Freedman tubing\/} $\hat F$
of $F$ is a surface obtained from $F$ by adding some annuli on $\bdd
M$ with boundary on $\bdd F$.  This idea has been used in several
important works, see [CLR, CL1, CL2, Li].  In particular, it was first
proved by Cooper and Long [CL2] that a Freedman tubing of an embedded,
geometrically finite surface is essential if the tubes are long
enough.  A combinatorial proof has been given by Li [Li], which also
yields an upper bound of tube length in terms of genus and number of
boundary components of $F$.  They have also shown that the tubed
surface survives most Dehn fillings, which, combined with a result of
Culler and Shalen [CS], implies that all but finitely many Dehn
filling spaces of a hyperbolic manifold contain an immersed surface.

Define the {\it wrapping number\/} $\text{wrap}(A)$ of an annulus $A$
on a torus $T$ to be the minimum algebraic intersection number between
$A$ and all points of $T$.  If $\hat F$ is a Freedman tubing of $F$,
define $w(\hat F, F)$ to be the minimum of $\text{wrap}(A_i)$ over all
components $A_i$ of $\hat F - F$.  The following theorem generalizes
the above result to immersed essential surfaces.

\proclaim{Theorem 5.7} Let $F$ be a geometrically finite surface in a
compact hyperbolic 3-manifold $W$.  Then there is a constant $K$ such
that if $\hat F$ is a Freedman tubing of $F$ with $\text{wrap}(\hat F,
F) \geq K$, then $\hat F$ is $\pi_1$-injective in $W$.
\endproclaim

The assumption that $F$ be geometrically finite is necessary,
otherwise $F$ would be a virtual fiber, and hence no Freedman tubing
of it would be essential.  Immersed surfaces are much more abundant
than embedded ones.  For example, Oertel [Oe] and Maher [Ma] showed
that in certain manifolds all slopes are realized as boundary slopes
of immersed essential surfaces, while Hatcher [Ha] showed that there
are only finitely many boundary slopes of embedded surfaces in these
manifolds.  More immersed surfaces can be obtained by projecting to
$M$ embedded surfaces in covering spaces of $M$.  The boundary of such
a surface may have several different slopes on the same torus
component of $\bdd M$.  Our theorem applies to such surfaces as well,
and there is no restriction on the orientability of $F$ or $\hat F$.
When $\bdd M$ is a set of tori and $F$ is a proper surface, a
essential Freedman tubing is automatically geometrically finite
because it has accidental parabolics, hence by Theorem 1.1 it
survivesill survive most Dehn fillings.

The idea of our proof is to use area estimation to show that certain
curves in a negatively curved space are nontrivial.  In section 2 we
will use some results in minimal surface theory to show that if a
piecewise geodesic curve $\a$ is trivial in a negatively curves
3-manifold $M$, then it bounds a disk whose intersection with the
hyperbolic part of $M$ has area bounded above by the total external
angle of $\a$.  In section 3 we give some estimation for areas of
surfaces in truncated hyperbolic cusps, using integral of certain
differential forms and Stokes theorem.  These result will then be used
in section 4 to show that curves in $M$ satisfying certain conditions
do not bound any disk, hence is nontrivial in $M$.  The essentiality
of surfaces in Dehn filling space and the essentiality of tubing
surfaces in hyperbolic manifolds follow from these results by showing
that all nontrivial curves on the surface satisfy those conditions.
In section 6 we will show that there is no universal upper bounds for
the bad fillings, and post several problems arisen in this research.

\medskip

{\it Definitions and conventions.}  All 3-manifolds in this paper are
assumed orientable.  Let $F: S \to M$ be a surface.  A point $p\in S$
is a {\it regular point\/} if $F$ is a local immersion at $p$.
Otherwise it is a {\it singular point.}  Almost all points of $S$ are
regular points since $F$ is assumed to be an immersion almost
everywhere.  We will use $F$ the same way as we would for embedded
surfaces.  Thus for example $\bdd F$ denotes the restriction of $F$ to
$\bdd S$, and if $N$ is a submanifold of $M$ then $F \cap N$ denotes
the restriction of $F$ on $F^{-1}(N) \subset S$, which is considered
as a subsurface of $F$ if $F^{-1}(N)$ is a subsurface of $S$.  By a
{\it curve\/} on a surface $F: S \to M$ we mean the composition
$F\circ \a$, where $\a: S^1 \to S$ is a closed curve on $S$.
Similarly if $\a: I \to S$ is an arc then $F\circ \a$ is called an
{\it arc\/} on $F$.  We say that the arc $F\circ \a$ has endpoints on
$\bdd F$ if $\bdd \a \subset \bdd S$, in which case it is a called a
proper arc.

Given two arcs or curves $\a, \b$ on $F$ or in a manifold $M$, we use
$\a \sim \b$ to denote that $\a, \b$ are homotopic.  Homotopy of arcs
and curves are different.  Two arcs $\a, \b$ are {\it homotopic\/} if
they are homotopic rel boundary in the usual sense, while two curves
are homotopic if they are freely homotopic.  A curve in a space is
{\it trivial\/} if it is null-homotopic.  An arc $\a$ on a surface $F$
is {\it essential\/} if it is not homotopic to an arc on $\bdd F$.

A surface $F: S \to M$ is {\it incompressible\/} if any nontrivial
curve on $F$ is also nontrivial in $M$.  Note that $F$ is
incompressible if and only if it is {\it $\pi_1$-injective\/}, that
is, $F_*: \pi_1 S_i \to \pi_1 M$ is an injective map for all
components $S_i$ of $S$.  A compact surface $F$ is {\it proper\/} if
$\bdd F \subset \bdd M$.  $F$ is {\it $\bdd$-incompressible\/} if no
essential arc of $F$ is homotopic in $M$ to an arc on $\bdd M$.  A
proper surface $F$ in $M$ is {\it essential\/} if it is
incompressible, $\bdd$-incompressible, and is not rel $\bdd F$
homotopic to a surface on $\bdd M$.

We refer the readers to [Th1] and [Mg] for basic concepts about
hyperbolic 3-manifolds.  In different sections below, $M$ may denote
either a complete hyperbolic manifold or a compact manifold with
interior a complete hyperbolic manifold.  If $M$ is a complete
hyperbolic manifold, the {\it injective radius\/} of a point $x$ in
$M$ is the supremum of radii of all embedded balls in $M$ centered at
$x$.  Denote by $M_{(0,\e]}$ the set of points which has injective
radius at most $\e$, and by $M_{[\e, \infty)}$ the set with injective
radius at least $\e$.  It is well known (see [Mg]) that when $\e$ is
sufficiently small, $M_{(0,\e]}$ is a set of cusps, in which case we
use $N = N_{\e}$ to denote the toroidal cusp components of
$M_{(0,\e]}$, and $T = T_{\e}$ the boundary tori of $N$.

The hyperbolic structure of $M$ induces a Euclidean metric on $T =
T_{\e}$.  If $\a$ is either a curve on $T$ or an arc in $M$ which is
homotopic to an arc on $T$, then $\a$ can be homotoped to a geodesic
$\a'$ on $T$.  Define $t(\a)$ to be the Euclidean length of $\a'$, and
call it the {\it $T$-length\/} of $\a$.  Notice that it depends only
on $\e$ and the homotopy class of $\a$.  If $\r$ is another curve or
arc on $T$, and $\r'$ the geodesic on $T$ homotopic to $\r$, then the
{\it $T$-length of $\a$ relative to $\r$\/}, denoted by $t_{\r}(\a)$,
is defined as
$$ t_{\r}(\a) = t(\a) |\sin \theta|, \tag{1-1} $$ where $\theta$ is
the angle between $\a'$ and $\r'$.  Geometrically, $t_{\r}(\a)$ is the
length of the orthogonal projection of $\a'$ to a line orthogonal to
$\r'$.  These notations will be used throughout the paper.

\medskip

{\it Acknowledgement.}  I would like to thank Charlie Frohman for some
useful conversations on this work, to him and Oguz Durumeric for helps
on minimal surface theory, and to Francis Bonahon and Darren Long for
some helpful communications.

\head {\S 2.  Minimal surfaces and the Plateau problem}
\endhead

Let $F: S \to M$ be a surface of type $S$ in a Riemannian manifold
$M$.  Recall that $F$ is assumed piecewise smooth.  In this
section we will also assume that $F$ is oriented.
If $\omega$ is a
differential 2-form of $M$, then by the restriction of $\omega$ to $F$
we mean the 2-form $F^*(\omega)$ on $S$ defined on all smooth points
of $F$, and the integral of $\omega$ on $F$ is defined as
$$ \int_F \omega = \int_S F^*(\omega).$$ Since $F$ is piecewise
smooth, this is well defined.

The Riemannian metric on $M$ induces a Riemannian metric on the set of
regular points of $F$, which determines a volume form $\omega_F$.
More explicitly, if $(u_1, u_2)$ is a local coordinate system of $S$
at a regular point $p$ of $F$ which is compatible with the orientation
of $S$, then the tangent vectors $\bdd_i = \frac{\bdd}{\bdd u_i} \in
T_pS$ are mapped to $F_*(\bdd _i)$ in $T_{F(p)} M$.  The Riemannian
metric of $M$ determines an inner product $\left< \cdot , \cdot
\right>$ on $T_{F(p)}M$.  Let $g_{ij} = \left< F_*(\bdd_i),
F_*(\bdd_j)\right>$.  Then
$$ \omega_F = \sqrt {\det(g_{ij})} \,\,\, du_1 \wedge du_2.$$
This is a well-defined 2-form on $S$.  Given a function $f(p)$ on
$S$, which we consider as a function on $F$, the integral of $f$
on $F$ is defined as
$$ \int _F f = \int _S f(p) \,\, \omega_F.$$
In particular, when $f = 1$, this defines the area of $F$:
$$ \area(F) = \int _F 1 = \int_S \omega_F.$$

If $M$ is of dimension two, then it has a volume form $\omega_M$, in
which case $\omega_F = \pm F^*(\omega_M)$, where the sign depends on
whether $F$ is orientation preserving or orientation reversing at that
point.  Given a 2-form $\omega$ on $S$ with local presentation $\omega
= \varphi \, du \wedge dv$, where $(u,v)$ is a local coordinate system
compatible with the orientation of $S$, we use $|\omega|$ to denote
the 2-form $|\varphi| du \wedge dv$.  Thus $\omega_F =
|F^*(\omega_M)|$ when $M$ is a surface.

We refer the readers to [Dc] for the definitions of curvatures and
second fundamental form of submanifolds.  Let $(h_{ij})$ be the second
fundamental form of $F$ at a regular point $p$, with respect to a
basis $(v_1, v_2)$ of $T_pF \subset T_p M$, then the Gauss formula
(cf.\ [Dc, p.130]) says
$$ K = \bar K(v_1, v_2) + \det (h_{ij}) = \bar K(v_1, v_2) +
h_{11}h_{22} - h_{12}^2$$ where $K$ is the curvature of $F$, and $\bar
K$ is the sectional curvature of $M$.  A continuous map $F: S \to M$
is a {\it minimal surface\/} if it is smooth in the interior of $S$,
and its mean curvature $h_{11} + h_{22}$ is always zero.  $F$ is not
required to be smooth on $\bdd S$.  Thus if $F$ is a minimal surface
then $h_{11}h_{22} \leq 0$, so from the above we have $K \leq \bar
K(v_1, v_2)$.  In particular, if $M$ is a hyperbolic manifold, which
by definition has constant sectional curvature $\bar K = -1$, then $ K
\leq -1$.

The classical Plateau problem asks if a Jordan curve in $\Bbb R^n$
bounds a surface of disk type with minimal area.  A solution to the
Plateau problem is necessarily a minimal surface, which is harmonic in
the interior of $D$, and is continuous on $D$.  The problem was first
solved by Douglas [Dg], and has been generalized by Morrey [Mr] to many
Riemannian manifolds.  The regularity of solutions has also been
deeply studied.  For our purpose, the following result suffices.

\proclaim{Lemma 2.1} Let $C$ be a null-homotopic, smooth, embedded
circle in a complete, negatively curved 3-manifold $M$ with hyperbolic
ends.  Then

(i) $C$ bounds a minimal surface $F: D^2 \to M$ of disk type, which
minimizes the area of all disk type surfaces bounded by $C$;

(ii)  $F$ is a smooth map on $D^2$;

(iii) if $K$ is the curvature function of $F$, and $\kappa$ the
geodesic curvature function of $C$ in $M$, then
$$\int _F K + \int _C \kappa \geq 2\pi.$$
\endproclaim

\proof (i) This follows from Morrey's solution of the Plateau problem
for Riemannian manifolds [Mr].  Morrey's result says that if $M$ is a
complete Riemannian manifold which is almost homogeneous, then any
null-homotopic curve $C$ in $M$ bounds a minimal surface which
minimizes the area of all disk type surfaces bounded by $C$.  Since we
have assumed that $M$ is complete and has hyperbolic ends, Morrey's
result applies.

(ii) This follows from Theorem 4 in Chapter 7 of [DHKW], which says
that the degree of smoothness of a minimal surface on its boundary $C$
is at least that of $C$ and $M$.  Since both $C$ and $M$ are assumed
smooth, the result follows.

(iii) We need the following Gauss-Bonnet theorem for minimal surfaces
with smooth boundary:
$$ \int_F K + \int_{\bdd F} \kappa_g = 2\pi + 2\pi \sum_{w \in
\sigma'} \nu(w) + \pi \sum_{w \in \sigma''} \nu(w),$$ where $\kappa_g$
is the geodesic curvature of $\bdd F$ on $F$, $\sigma'$ and $\sigma''$
the set of interior and boundary branch points, respectively, and
$\nu(w)$ the branch index of $w$, which is nonnegative.  For minimal
surfaces in $\Bbb R^n$, this is Theorem 1 in Chapter 7 of [DHKW], and
for minimal surfaces in Riemannian manifolds it is proved by Kaul [K].
The proof would be easy if one knows the smoothness of $F$ on its
boundary, which is (ii) above, and the local behavior of $F$ near its
branch points, which was done by Heinz and Hildebrandt [HH].  If $p$
is a regular point of $F$ on $\bdd F$, $\bold n$ the principal normal
vector of $\bdd F$ at $p$, $\bold n'$ the inward normal vector of
$p$ on $F$, and $\theta$ the angle between $\bold n$ and $\bold n'$,
then $\kappa_g = (\cos \theta) \, \kappa$.  Therefore we have
$\kappa_g \leq \kappa$, and the result follows.
\endproof

If $C$ is a piecewise geodesic curve, and $p$ is a corner point of
$C$, then by going around $C$ in a certain direction, we get two
tangent vectors $v_1, v_2$ at $p$.  The {\it external angle\/} of $C$
at $p$ is the angle between $v_1$ and $v_2$.  The {\it total external
angle\/} of $C$ is the sum of external angles at all the corner points
of $C$.

\proclaim{Proposition 2.2} Let $M$ be a complete negatively curved
3-manifold with hyperbolic ends, and let $M_h$ be a hyperbolic
submanifold of $M$.  Suppose $C$ is a piecewise geodesic in $M$ such
that $M$ is hyperbolic near all corners of $C$.  Let $\Theta$ be the
total external angle of $C$.  Then $C$ bounds a surface $F$ of disk type
in $M$ such that
$$\area(F \cap M_h) \leq \Theta - 2\pi.$$
\endproclaim

\proof At each corner $p$, let $D_p = \exp D_{\delta}$, where
$D_{\delta}$ is a disk of radius $\delta$ on the plane in $T_pM$
containing the two tangent vectors of $C$ at $p$, and $\exp$ the
exponential map.  Since $M$ is hyperbolic near $p$, by choosing
$\delta$ small enough we may assume that $D_p$ is an embedded totally
geodesic disk in $M$.  Let $\a'_1, \a'_2$ be the two geodesic segment
of $C \cap D_p$.  Choose a point $p_i$ in the interior of each
$\a'_i$, and let $\a_i$ be the subarc of $\a'_i$ connecting $p_i$ to
$p$.  Connect $p_1$ to $p_2$ by a smooth arc $\r_p$ such that
$(C-\a_1\cup \a_2) \cup \r_p$ is smooth in $D_p$, and $\r_p$ is
concave on the region $\Delta_p$ bounded by $\a = \a_1\cup \a_2\cup
\r_p$.  See Figure 2.1.

\bigskip
\leavevmode

\centerline{\epsfbox{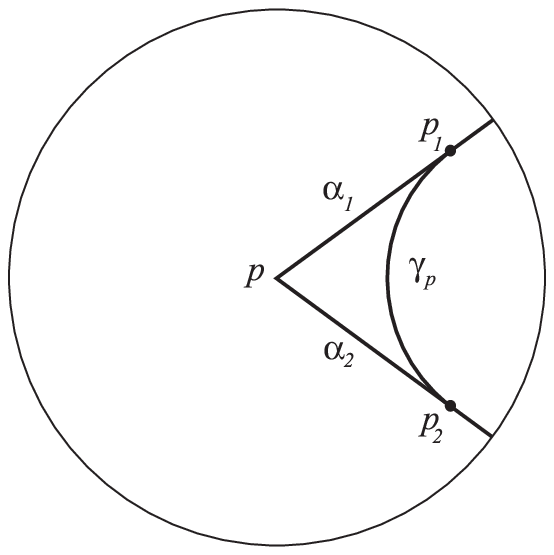}}
\bigskip
\centerline{Figure 2.1}
\bigskip

Since $D_p$ is totally geodesic, the curvature $\kappa$ of $\r_p$ in
$M$ is the same as that in $D_p$.  Since it is concave as a boundary
curve of $\Delta_p$, its curvature $\kappa_g$ as boundary curve of
$\Delta_p$ is $- \kappa$.  The total external angle of $\bdd \Delta_p$ is
$2\pi + \theta(p)$, where $\theta(p)$ is the external angle of $C$ at
$p$.  Therefore by the Gauss-Bonnet theorem applied to $\Delta_p$, we
have
$$\int_{\Delta_p} (-1) + \int_{\bdd \Delta_p} \kappa_g + (2\pi +
\theta(p)) = 2\pi.$$ The first integral is $-\area(\Delta_p)$, and the
second equals $-\int_{\r_p} \kappa$.  Hence
$$\area(\Delta_p) + \int_{\r_p} \kappa = \theta(p).$$ Let $C'$ be the
smooth curve obtained from $C$ by replacing $\a_1\cup \a_2$ with
$\r_p$ at each corner $p$, and let $F'$ be the minimal surface bounded
by $C'$ as given in Lemma 2.1.  Then $F = F' \cup (\cup \Delta_p)$ is
a surface bounded by $C$.  Since the curvature $K$ of $F'$ satisfies
$K\leq -1$ in $M_h$ and $K<0$ elsewhere, by Lemma 2.1(3) and the above
we have
$$ \align
\area(F \cup M_h) &\leq \sum_p \area(\Delta_p) + \area(F' \cap M_h)
\leq \sum_p \area(\Delta_p) - \int_{F'} K \\
& \leq \sum_p \area(\Delta_p) + \int_{\bdd F'} \kappa \, - 2\pi
= \sum_p\left[\area(\Delta_p) + \int_{\r_p} \kappa \right] - 2\pi  \\
& = \sum \theta(p) - 2\pi = \Theta - 2\pi.
\endalign$$
\endproof

\remark {Remark 2.3} (1) Charles Frohman pointed out that when $M$ is
hyperbolic, Proposition 2.2 can be proved easily by considering a disk
bounded by $C$ which is a union of totally geodesic triangles.  Thus
the above proof using minimal surface theory is necessary only if $M$
is negatively curved but not hyperbolic.

(2) Proposition 2.2 would follow more directly if we had a
Gauss-Bonnet type formula for minimal surfaces with boundary a smooth
curve with corners.  It should look like:
$$ \int_F K + \int_{\bdd F} \kappa_g + \sum_{p\in \sigma'''} (\pm
\theta_p) = 2\pi + 2\pi \sum_{w \in \sigma'} \nu(w) + \pi \sum_{w \in
\sigma''} \nu(w)$$ where $\sigma'''$ is the set of corner points, and
$\theta_p$ the external angle of $C$ at $p$.  Note that negative sign
could appear before $\theta_p$ if $p$ is a branch point.  The formula
could be proved in the usual way if we know the local behavior of $F$
near the corner points, which was done in Chapter 8 of [DHKW] in the
special case that $F$ is in Euclidean space.  Unfortunately I cannot
find a reference for either the formula or the local behavior near
corners of a minimal surface $F$ in a Riemannian manifold.  \endremark

\head 3. Area estimation for surfaces in truncated hyperbolic cusps
\endhead

Throughout this paper, we will always consider the hyperbolic space
$\H^3$ as in the upper half space model.  Denote by $\H^3_1$ the
hyperbolic horoball $\{(x,y,z) \, |\, z \geq 1\}$.  For $b>1$, denote
by $\H^3_{1,b}$ the subset of $\H^3_1$ where $z \leq b$.  Consider
$\H^2$ as the subset of $\H^3$ corresponding to the $yz$-plane.
Define $\H^2_1 = \H^3_1 \cap \H^2$, and $\H^2_{1,b} = \H^3_{1,b} \cap
\H^2$.  For simplicity, we use $(y,z)$ to denote a point $(0,y,z)$ in
$\H^2$.

Consider the following subset $R_1(a,b)$ and $R_2(a,b)$ of $\H^2_1$ as
shown in Figure 2.1, where $R_1(a,b)$ is a Euclidean rectangle, and
$R_2(a,b)$ is the intersection with $\H^2_{1,b}$ of a Euclidean disk
which is centered at the origin and intersects the horizontal line at
$z=1$ in an arc of length $a$.  Thus it has radius $\sqrt{1 +
(a/2)^2}$.  More explicitly, we have
$$ \align R_1(a,b) & = \{(y,z) \in \H^2 \,\,|\,\, 0\leq y \leq a,
\,\,1\leq z \leq b\},\\ R_2(a,b) & = \{(y,z) \in \H^2 \,\,|\,\,
1\leq z \leq b, \,\, y^2 + z^2
\leq 1 + (\frac a2)^2 \}.  \endalign
$$
Define a function $\eta(x)$ by
$$ \eta(x) = x - 2 \arctan{\frac x2}.$$

\bigskip
\leavevmode

\centerline{\epsfbox{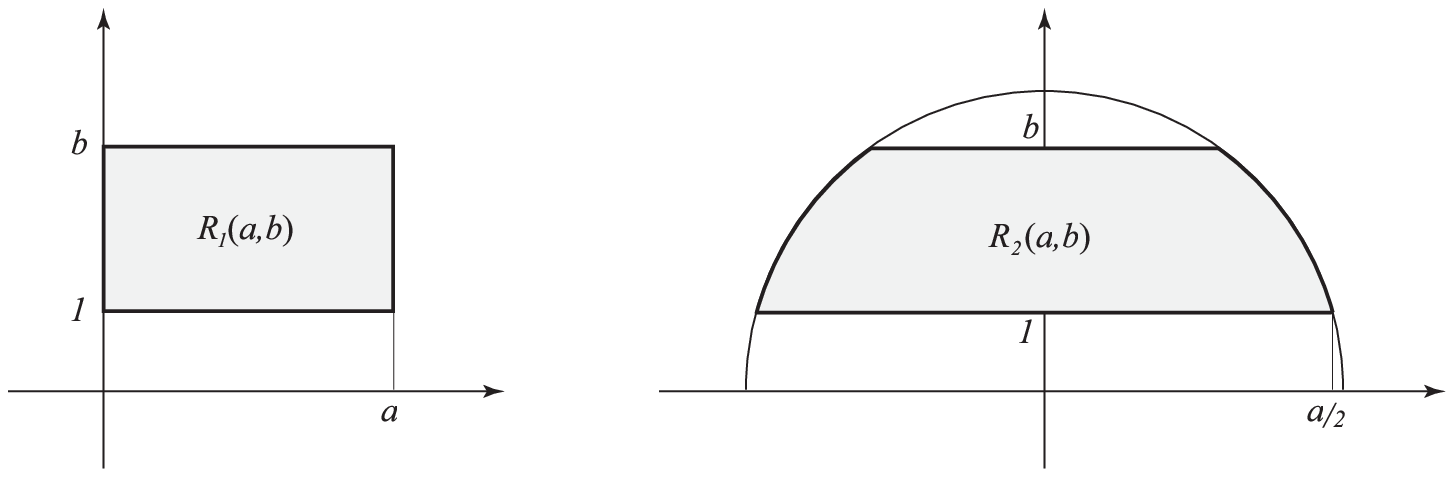}}
\bigskip
\centerline{Figure 3.1}
\bigskip

\proclaim{Lemma 3.1} (1) $\area (R_1(a,b)) = a(1 - 1/b)$.

\medskip

(2) $\displaystyle \area(R_2(a,b)) = \left\{
\matrix \format\l&\qquad\l\\
\eta(a)  &  b^2 \geq 1 + \frac {a^2}4 \\
\eta(a) - \eta(\frac{\sqrt{1 + a^2/4 - b^2}}b)
& b^2 \leq 1 + \frac {a^2}4
\endmatrix
\right.
$

\medskip

(3)  If $a\geq 3\pi$ and $b\geq 5$, then $\area(R_2(a,b)) > 2\pi$.
\endproclaim

\proof These would follow from the Gauss-Bonnet theorem and the fact
that a horizontal line in $\H^2$ at height $b$ has curvature $1/b$,
with normal vector pointing upward.  The following is a direct
calculation.

(1)  $\displaystyle \area(R_1(a,b)) = \int_0^a dy \int_1^b
\frac 1 {z^2} \, dz = a(1-\frac 1b)$.

(2)  Let $r = \sqrt{1+ \frac {a^2}2}$.  First assume $b\geq r$.  Then
$$ \align \area(R_2(a,b)) &= \iint_{R_2(a,b)} \frac 1{z^2} \, dy\, dz
= 2 \int_1^r \frac 1{z^2} \, dz \int_0^{\sqrt{r^2 - z^2}} dy \\ & = 2
\int_1^r \frac{\sqrt{r^2 - z^2}}{z^2} \, dz = 2\left[ \arctan
\frac{\sqrt{r^2 - z^2}}z - \frac{\sqrt{r^2 - z^2}}z \right]_1^r \\ &=
2(\frac a2 - \arctan \frac a2) = \eta(a).  \endalign
$$ When $b < r$, $R_2(a,b) = R_2(a, \infty) - R'$, where $R'$ is the
subregion of $R_2(a,\infty)$ above the line $z=b$.  The transformation
$(y,z) \to (y/b, z/b)$ is a hyperbolic isometry, which maps $R'$ to
the region $R_2(\frac {\sqrt{r^2 - b^2}}b, \infty)$, so the result
follows from the above.

(3) From the definition it is clear that $\area(R_2(a,b))$ is an
increasing function of both $a$ and $b$.  Since $a\geq 3\pi$ and $b^2
\geq 5^2 > 1 + (1.5 \pi)^2 = 1+a^2/4$, by (2) we have
$$ \area(R_2(a,b)) \geq \area(R_2(3\pi, b)) = 3\pi - 2 \arctan \frac
{3\pi}2 > 2\pi.$$
\endproof

The hyperbolic metric on $\H^3$ induces a Euclidean metric on the
Euclidean plane $P = \bdd \H^3_1$.  Recall from [Th1] that a
hyperbolic cusp $N$ of toroidal type is isometric to $\H^3_1 / G$ for
some Euclidean translation group $G$ of $P$ of rank $2$.  Denote by
$T$ the boundary torus of $N$, and by $N^b$ the truncated cusp
$\H^3_{1,b} /G$.  We allow $b= \infty$, in which case $N^b = N$.

If $\r$ is a nontrivial closed curve on $T$, then there is a totally
geodesic annulus $A_{\r}$ in $N^b$ perpendicular to the boundary, such
that $A_{\r} \cap T$ is homotopic to $\r$.  More explicitly, up to
rotation and translation of $\H^3_1$ we may assume that $\r$ lifts to
an arc on $P = \bdd \H^3_1$ with both endpoints on the $y$-axis.  Let
$A'_{\r}$ be the annulus obtained from $R_1(t(\r), b)$ by identifying
the two vertical lines.  Then the quotient map $q$ from $\H^3_{1,b}$ to
$N^b$ induces a map on $A'_{\r}$, which we define as the surface
$A_{\r}$ in $N^b$.  Notice that if $\r = k \b$ in $H_1(T)$, then
$A_{\r}$ is a $k$-fold cover of $A_{\b}$.  By Lemma 3.1(1) we have
$$ \area (A_{\r}) = \area (R_1(t(\r), b) = t(\r)(1 - \frac 1b).$$

Let $F$ be a surface of type $S$ in $N$ or $\H^3$.  We would like to
estimate the area of $F$.  Consider the 2-form
$$ \omega = \frac 1{z^2} \,\, dy\wedge dz$$ on $\H^3$.  Notice that
its restriction to $\H^2$ is the standard volume form $\omega_{\H^2}$
of $\H^2$, and if we denote by $p: \H^3 \to \H^2$ the Euclidean
orthogonal projection $p(x,y,z) = (y,z)$, then $\omega =
p^*(\omega_{\H^2})$.  Therefore if $F: S \to \H^3$ is a surface in
$\H^3$ then
$$ \area(p\circ F) = \int_S |(p\circ F)^*(\omega_{\H^2})| =
\int_S |F^*(\omega)| = \int_F |\omega| \geq |\int_F \omega\, |$$

The map $p$ is area non-increasing, so the area of $F$ is at least
that of $p\circ F$.  In fact, more is true.  Recall that $\omega_F$
denote the volume form of $F$ induced by the Riemannian metric of
$\H^3$.

\proclaim{Lemma 3.2} Let $F: S \to \H^3$ be a surface in $\H^3$.  Let
$\theta(p)$ be the angle between the normal vector of $F$ at a regular
point $p$ and the positive $x$-axis.  Then $ F^*(\omega) = \cos
\theta(p)\,\, \omega_F.$ In particular, if $F$ is a Euclidean planar
surface in $\H^3$ (so $\theta$ is a constant), then $$\int_F \omega =
(\cos \theta) \area(F).$$  \endproclaim

\proof Let $(u,v)$ be a local coordinate system at a regular point
$p$.  Then $\bold n = F_u \times F_v$ is a normal vector of $T_p F$.
Put $F = (x(u,v), y(u,v), z(u,v))$, and $\bold n = n_1 \bold i + n_2
\bold j + n_3 \bold k$.  Then $n_1 = y_u z_v - y_v z_u$ and $\cos
\theta(p) = n_1 / || \bold n ||$.  Use $\bold a \cdot \bold b$ to
denote the dot product of two vectors in $\Bbb R^3$.  Then $g_{11} =
\left< F_u, F_u \right> = \frac 1{z^2} F_u \cdot F_v$.  Similarly for
the other $g_{ij}$.  Thus
$$\det (g_{ij}) = \frac 1 {z^4} [ (F_u \cdot F_u) (F_v \cdot F_v)
- (F_u \cdot F_v)^2 ] = \frac 1 {z^4} ||F_u \times F_v ||^2.$$
Hence
$$ \omega_F = \sqrt{\det (g_{ij})} \, du \wedge dv = \frac 1 {z^2}
|| F_u \times F_v || \, du \wedge dv = \frac {||\bold n||}
{z^2} \, du \wedge dv.$$
On the other hand, we have
$$ \align F^*(\omega) &= F^*(\frac 1 {z^2} \, dy \wedge dz) = \frac 1{z^2}
(y_u du + y_v dv) \wedge (z_u du + z_v dv) \\
& = \frac 1{z^2}
(y_u z_v - y_v z_u) \, du \wedge dv = \frac {n_1}{z^2} \, du \wedge
dv = \cos \theta(p) \, \omega_F. \endalign $$
\endproof

Now consider a hyperbolic cusp $N = \H^3_1/G$, with torus boundary $T
= \bdd N$.  Clearly $\omega$ is invariant under Euclidean
translations, hence it induces a 2-form $\omega_N$ on $N$.  Suppose
$F: S \to N$ is a surface in $N$.  Since Lemma 3.2 is a local property,
we still have
$$ F^*(\omega_N) = \cos \theta(p)\,\, \omega_F$$
where $\theta(p)$ is the angle between the normal vector of $F$ at $p$
and a vector in $T_pN$ whose lifting to $\H^3_1$ points to the positive
$x$-axis direction.

\proclaim{Lemma 3.3} Let $F, F_1, F_2: S \to N^b$ be compact, oriented
surfaces in $N^b$ with boundary on $\bdd N^b$.  Then

(1) $\displaystyle \area(F) = \int _S \omega_F \geq \int_F |\omega_N|
\geq |\int_F \omega_N|$;

(2) If $[F_1] = [F_2] \in H_2(N^b, \bdd N^b)$, then $\displaystyle
\int _{F_1} \omega = \int_{F_2} \omega$;

(3) If $[\bdd F \cap T] = [\r] \neq 0 \in H_1(T)$ and $\delta$ is a
geodesic arc on $T$ which lifts to an arc on $\bdd \H^3_1$ parallel to
the $x$-axis, then $\displaystyle |\int_F \omega_N | =
t_{\delta}(\r)(1 - \frac 1b)$.

(4) If $[\bdd F \cap T] = [\r] \neq 0 \in H_1(T)$, then $\displaystyle
\area(F) \geq t(\gamma) (1-\frac 1b)$.  \endproclaim

\proof (1) We have
$$ \int_S \omega_F = \int_S |\omega_F| \geq \int_S |\cos \theta(p) \,
\omega_F| = \int_S |F^*(\omega_N)| = \int_F |\omega_N|.$$

(2) Notice that
$$d(\omega) = d\,(\frac 1{z^2} \,dy \wedge dz) = \frac {-2}{z^3} \,dz
\wedge dy \wedge dz = 0,$$ so $\omega$ is a closed form.  Since
$\omega_N$ is induced from $\omega$, it is also a closed form.  Denote
by $\bar F_2$ the surface $F_2$ with orientation reversed.  The
assumption means that there is a surface $F_3$ on $\bdd N^b$, such
that $\hat F = F_1\cup \bar F_2 \cup F_3$ is a closed oriented surface
which is null homologous in $N^b$.  Therefore there is an oriented
3-manifold $W$ and a map $f: W \to N^b$ with $f|_{\bdd W} = \hat F$.
By Stokes theorem, we have
$$ \align \int_{\hat F} \omega_N &= \int_{\bdd W} f^*(\omega_N) =
\int_W d(f^*(\omega_N)) \\ &= \int_W f^*(d\, \omega_N)) = \int_W 0 =
0.\endalign$$ Since $F_3$ lifts to a horizontal planar surface in
$\H^3$, by Lemma 3.2 $$\int _{F_3} \omega_N = (\cos \frac{\pi}2)\,
\area(F_3) = 0.$$ Therefore
$$ 0 = \int_F \omega_N = \int_{F_1} \omega_N + \int_{\bar F_2} \omega_N
+ \int_{F_3} \omega_N = \int_{F_1} \omega_N - \int_{F_2} \omega_N$$
and the result follows.

(3) $F$ is homologous to the surface $A_{\r}$ defined above, which
lifts to a region on a vertical plane in $\H^3_1$.  Let $\theta'$ be
the angle between $\d$ and $\r$.  Then the acute angle $\theta$
between the normal vector of $A_{\r}$ and the $x$-axis satisfies
$\theta = |(\pi/2) - \theta'|$.  Hence by the definition of
$t_{\d}(\r)$ is section 1, we have $\cos \theta = |\sin \theta'| =
t_{\delta}(\r) / t(\r)$.  It follows from (2) and Lemma 3.2 that
$$|\int_F \omega_N| = |\int_{A_{\r}} \omega_N| = (\cos \theta) \,
\area(A_{\r}) = (\cos \theta)\, t(\r)(1 - \frac 1b)
= t_{\delta}(\r) (1 - \frac 1b).$$

(4) Choose a coordinate system of $\H^3$ so that the geodesic $\r'$
homotopic to $\r$ lifts to the $y$-axis.  Let $\delta$ be an arc
perpendicular to $\r'$.  Then
$$ \area(F) \geq |\int_F \omega_N| = t_{\delta}(\r) (1 - \frac 1b) =
t(\r) (1 - \frac 1b).$$ \endproof

\proclaim{Lemma 3.4} Let $\b$ be an arc on $T$, and let $\a$ be the
geodesic segment in $N$ homotopic to $\b$.  Let $F$ be a compact,
oriented surface in $N^b$ such that $\bdd F \cap \Int N^b = \a \cap
\Int N^b$.  Put $[\r] = [(\bdd F\cap T) \cup \b] \in H_1(T)$.

(1) If $[\r] = 0$, then $\area(F) \geq \area(R_2(t(\a), b)))$.

(2) If $[\r] \neq 0$, then $\area(F) \geq t_{\d}(\r)(1- \frac 1b) -
t_{\d}(\a)$ for all slopes $\d$ on $T$.  \endproclaim

\proof Without loss of generality we may assume that $\b$ is a
geodesic on $T$.  Choose a coordinate system of $\H^3$ so that the
lifting of $\b$ is an arc on the $y$-axis, and is symmetric about the
$z$-axis.  Let $R$ be the image of $R_2(t(\a), b)$ in $N^b$ under the
projection map.  Then $\hat F = F \cup R$ is a properly embedded
surface in $N^b$.  Notice that $(\bdd \hat F) \cap T = (\bdd F \cap
T)\cup \b$, so it is homologous to $\r$.

In case (1), the surface $\hat F$ is null-homologous in $H_2(N^b, \bdd
N^b) \cong H_1(T)$, hence by Lemma 3.3(2) we have $\int_{\hat F}
\omega_N = 0$.  Thus $\int_F \omega_N = -\int_R \omega_N$.  Since $R$
lifts to the region $R_2(t(\a), b)$ on $\H^2$, we have
$$\area(F) \geq |\int_F\omega_N| = |\int_R \omega_N | =
\area(R_2(t(\a), b)).$$

In case (2), rechoose the coordinate system so that the geodesic on
$T$ homotopic to $\d$ lifts to the $x$-axis.  Then $R$ lifts to a
surface $\tilde R$ which is a rotation of $R_2(t(\a), b)$ by an angle
$\theta$.  As in the proof of Lemma 3.3(3), we have $\cos \theta =
t_{\d}(\a)/t(\a)$, so by Lemma 3.2 we have
$$\int_R \omega_N = (\cos \theta)\, \area(R) = \frac{t_{\d}(\a)}{t(\a)}
\area(R_2(t(\a), b)).$$
By Lemma 3.3(3), $ |\int_{\hat F} \omega_N| = t_{\d}(\r)(1- 1/b)$.
Therefore,
$$\align
\area(F) &\geq |\int_F \omega_N |
\geq |\int_{\hat F} \omega_N - \int_R \omega_N |
\geq |\int_{\hat F} \omega_N|  - |\int_R \omega_N | \\
& \geq t_{\d}(\r)(1 - \frac 1{b}) - \frac{t_{\d}(\a)}{t(\a)}
\area(R_2(t(\a), b)) > t_{\d}(\r)(1 - \frac 1{b}) - t_{\d}(\a).
\endalign
$$
The last inequality is because by Lemma 3.1(2) we have $\area(R_2(t(\a),
b)) < t(\a)$.
\endproof

\head {\S 4. Nontrivial curves in negatively curved manifolds}
\endhead

Let $M$ be a complete hyperbolic manifold of finite volume.  Let $N$
be a set of mutually disjoint cusps of $M$.  Let $M_0 = M - \Int N$.
Put $T = \bdd M_0 = \bdd N$, which is a union of tori.  We may choose
$N$ so that $T$ lifts to a set of horospheres in $\H^3$, hence it has
a Euclidean metric induced by the hyperbolic metric of $M$.

A geodesic arc $\a$ in $M$ with endpoints on $T$ is said to be of {\it
type I\/} if a neighborhood of $\bdd \a$ lies in $M_0$, and of {\it
type II\/} if $\a \subset N$.  Notice that a geodesic arc may be
neither of type I nor of type II, but we will not consider such arcs.

\proclaim{Theorem 4.1} Let $M$ be a complete hyperbolic 3-manifold,
and let $M_0, N$ be as above.  If $\a = \a_1 \cup ... \cup \a_{2p}$ is
a closed curve such that (i) each $\a_{2i+1}$ is a geodesic arc of
type I, (ii) each $\a_{2i}$ is a proper arc in $N$, and (iii)
$t(\a_{2i}) \geq 2\pi$ for $i<p$, then $\a$ is nontrivial in $M$.

Moreover, if each $\a_{2i-1}$ has both endpoints perpendicular to $T$,
then (iii) can be replaced by (iii') $t(\a_{2i}) \geq \pi$ for $i<p$.
\endproclaim

\proof If the theorem were not true, we can choose a curve $\a$ as in
the theorem, such that $\a$ is null-homotopic in $M$, and $p$ is
minimal among all such curves.  By a homotopy we may assume
that all $\a_{2i}$ are geodesics in $N$, so they are type II arcs.
Now $\a$ is a piecewise geodesic curve with $2p$ corners, so its total
external angle is less than $2p\pi$.  By Proposition 2.2, it bounds a
surface $F: D^2 \to M$ of disk type, such that $\area(F) <
\Theta - 2\pi$, where $\Theta$ is the total external angle of $\a$.

By a small perturbation we may assume that $F$ is transverse to $T$.
Then $A = F \cap T$ is a compact 1-manifold in $F$.  Recall our
convention that we will treat $F$ the same way as an embedded surface.
Thus for example $A$ is really the restriction of $F$ on the
1-manifold $F^{-1}(T)$ in $D$, and by a disk cut off by a component of
$A$ we really mean the restriction of $F$ to a disk in $D$ cut off by
the corresponding component of $F^{-1}(T)$.

We claim that each arc component of $A$ is outmost in the sense that
it cuts off a disk $\Delta$ on $F$ containing no other {\it arc\/}
components of $A$.  (Note that $\Delta$ could contain some circle
components of $A$.)  Assuming otherwise, let $c$ be a component which
is not outmost.  Now $c$ is an arc on $T$, whose boundary cuts $\a$
into two arcs $\a'$ and $\a''$.  One of $\a'$ and $\a''$, say $\a'$,
does not contain $\a_{2p}$, so
$\a'\cup c$ satisfies the condition of the theorem with smaller
$p$, and is null-homotopic in $M$ because it bounds a subdisk of $F$.
This contradicts the minimality of $p$, completing the proof of the
claim.

Now let $\Delta$ be an outmost disk cut off by an arc component $c$ of
$A$.  Then $\Delta \cap \bdd D$ is one of the arcs $\a_i$ in $\a$.  We
have assumed above that $\a_i$ is a geodesic, so $\a_i$ being
homotopic to the arc $c$ on $T$ implies that $\a_i$ is in $N$, that
is, $i$ is an even number.  Hence we can label the outmost disks as
$\Delta_1, ..., \Delta_p$, with $\Delta_i\cap \bdd D = \a_{2i}$.

Recall that $\Delta_i$ may contain some circle components of $A$.  Let
$Q$ be the component of $\Delta_i$ cut along $A$ which contains $\bdd
\Delta_i$.  Since $M$ is hyperbolic, $N$ is $\pi_1$-injective in $M$,
hence each boundary component of $Q$ is null-homotopic in $N$ because
it bounds a disk in $M$.  Let $\b$ be an arc on $T$ homotopic to
$\a_{2i}$.  Then $(\bdd Q \cap T) \cup \b = (\bdd Q - \a_{2i}) \cup
\b$ is null-homologous on $T$ because $\b$ is homotopic to the arc
component of $\bdd Q \cap T$ and the circle components of $\bdd Q \cap
T$ are also null-homotopic.  Therefore, by Lemma 3.4(1) (with $b =
\infty$) and Lemma 3.1(2), for each $i<p$ we have
$$\area(\Delta_i) \geq \area(Q) \geq \area[R_2(t(\a_{2i}),\infty)] =
t(\a_{2i}) - 2 \arctan \frac {t(\a_{2i})}2$$ for all $i< p$.  Notice
that $\arctan (t(\a_{2i})/2)$ is the angle between $\a_{2i}$ and $T$.
Denote by $\theta_i$ the external angle at the corner between $\a_i$
and $\a_{i+1}$.  Since $a_j$ are of type I for odd $j$, we have
$\theta_{2i-1}, \theta_{2i} \leq \pi - \arctan (t(\a_{2i})/2)$, so the
above inequality together with the assumption $t(\a_{2i}) \geq 2\pi$
implies that $\area(\Delta_i) \geq \theta_{2i-1} + \theta_{2i}$ for
$i<p$.  Therefore
$$\area(F) > \sum_{i=1}^{p-1} \area(\Delta_i) \geq \sum_{j=1}^{2p-2}
\theta_j > \Theta - 2\pi.$$ Since $F$ is chosen to have area less than
$\Theta - 2\pi$, this is a contradiction.

If all $\a_{2i-1}$ have endpoints perpendicular to $T$, then $\theta_i
\leq \frac{\pi}2 - \arctan{(t(\a_{2i})/2)}$, so the assumption
$t(\a_{2i}) \geq \pi$ suffices to lead to a contradiction.
\endproof

We now consider Dehn fillings on $M$.  Recall that $N$ is a set of
disjoint cusps, and $M_0 = M - \Int N$.

Assume $t(\r_i) > 2\pi + 1$ for each $i$.  Choose $b_i$ so that the
geodesic curve $\r'_i$ on $T'_i = \bdd N^{b_i}_i - T_i$ isotopic to
$\r_i$ in $N_i$ has length $2\pi+1$.  Choose a coordinate for $\H^3$
so that the geodesic on $T_i$ homotopic to $\r_i$ lifts to the
$y$-axis.  Then the upper edge of $R_1(t(\r_i), b_i)$ is projected to
$\r'_i$.  Since the upper edge has hyperbolic length $t(\r_i)/b_i$, we
have
$$ b_i = \frac {t(\r_i)}{2\pi + 1}. \tag{4-1} $$

Denote by $N_i(\r_i)$ the manifold obtained by gluing a solid torus
$V_i$ to $N_i^{b_i}$ along $T'_i$ so that $\r'_i$ bounds a meridian
disk in $V_i$.  Put $N^b = \cup N_i^{b_i}$, $N(\r) = \cup N_i(\r_i)$,
and $M(\r) = M_0 \cup N(\r)$.  The manifold $M(\r)$ is the Dehn
filling space of $M$ (or more precisely, of $M_0$) along the multiple
slope $\r$.  By the $2\pi$-theorem of Gromov-Thurston [GT], $M(\r)$
has a negatively curved metric which coincides with the original
hyperbolic metric in a neighborhood of $M_0 \cup N^b$.  We will assume
below that $M(\r)$, $N(\r)$ and $V = \cup V_i$ are endowed with such a
metric.  Let $C_i$ be the core of $V_i$.  The identity map on $M_0
\cup N^b$ extends to a homeomorphism $M \cong M(\r) - \cup C_i$.  We
will always (topologically) identify $M$ with $M(\r)-\cup C_i$ in this
way; in particular, each curve $\a$ in $M$ is also a curve in $M(\r)$.

\proclaim{Lemma 4.2}  Let $K > 2\pi+1$ be a constant, and let $\r =
(\r_1, ..., \r_n)$ be a multiple slope on $T$ such that $t(\r_i) \geq
K$ for all $i$.  Let $D$ be a surface of disk type in $M(\r)$ such
that $\bdd D \subset T$, and $D$ is transverse to $T$.  If $\bdd D$ is
nontrivial on $T$, then $\area(D \cap N^b) \geq K-(2\pi+1)$.
\endproclaim

\proof Let $Q$ be the component of $D$ cut along $T$ containing $\bdd
D$.  If some component of $\bdd Q - \bdd D$ is nontrivial in $T$, then
by induction the subdisk $D'$ of $D$ bounded by this curve has
$\area(D' \cap N^b) \geq K - (2\pi+1)$, and we are done.  So assume
that all components of $\bdd Q - \bdd D$ are trivial on $T$.  If $Q$
were in $M_0$ then the above would imply that $\bdd D$ is
null-homotopic in $M_0$, contradicting the incompressibility of $T$ in
$M_0$.  Therefore $Q$ is contained in $N_i(\r_i)$ for some component
$N_i$ of $N$.  The above assumption means that each component of $\bdd
Q - \bdd D$ bounds a disk on $T_i$, hence $\bdd D$ is null-homotopic
in $N_i(\r_i)$.  Thus $[\bdd Q \cap T_i] = [\bdd D] = [k\r_i] \in
H_1(T_i)$, and $k\neq 0$ because $\bdd D$ is assumed nontrivial on
$T_i$.  Hence by Lemma 3.3(4) we have
$$ \area(D\cap N^b) \geq t(k\r_i)(1 - \frac 1{b_i}) = |k|(t(\r_i)
-\frac {t(\r_i)}{b_i}) \geq t(\r_i) - (2\pi+1) .$$ The last inequality
follows because $k \neq 0$, and because by (4-1) we have $t(\r_i) =
(2\pi+1)b_i$.  \endproof

\proclaim{Theorem 4.3} Let $\r = (\r_1, ..., \r_n)$ be a multiple
slope on $T$ such that $t(\r_i) \geq 12\pi$ for all $i$.  Let $\a =
\a' \cup \a''$ be a curve in $M$ such that either $\a''$ is a closed
geodesic and $\a' = \emptyset$, or $\a''$ is a type I geodesic arc and
$\a'$ is an arc in $N$.  If each component $\b$ of $\a'' \cap N_i$
satisfies $t_{\d}(\b) \leq t_{\d}(\r_i) - 5\pi$ for some slope $\d$ on
$T$, then $\a$ is nontrivial in $M(\r)$.  \endproclaim

\proof If $\a$ is a geodesic in $M_0$, then it remains a geodesic in
the negatively curved manifold $M(\r)$, hence is nontrivial.  (This is
well known, and also follows from Lemma 2.1(3) because $K <0$ and
$\kappa=0$.)  Therefore by choosing a component of $\a'' \cap N$ as
$\a'$ if necessary, we may always assume that $\a''$ is a type I
geodesic.  Put $\a'' = \a_1 \cup ... \cup \a_{2p-1}$.  Then $\a_{2j}$
lie in $N$, and $\a_{2j-1}$ are in $M_0$.  Assume the result is false, and let $\a$ be as in the theorem so that $\a$ is null
homotopic in $M(\r)$, and $p$ is minimal among all such curves.

Modify $\a$ as follows.  For each $\a_{2i}$ which has nontrivial
intersection with the Dehn filling solid tori $V_i$, homotope
$\a_{2i}\cap V_i$ to a geodesic segment $\a'_{2i}$ in $V_i$, and
denote the resulting arc $(\a_{2i} \cap N^b) \cup \a'_{2i}$ by
$\b_{2i}$.  Since $b_i = t(\r_i)/(2\pi +1) \geq 12\pi/(2\pi+1) > 5$,
from Figure 3.1 we see that such modification happens only if
$$ t(\a_{2i}) > 2 \sqrt{b_i^2 - 1} > 2 \sqrt{24} > 3\pi.$$ Let $r$ be
the number of arcs which have been modified.  Next, deform $\a'$ to a
geodesic $\b_{2p}$ in $N(\r)$.  For simplicity, write $\b_i = \a_i$
for the other arcs.  The curve $\b = \b_1 \cup ... \cup \b_{2p}$ is
now a piecewise geodesic in $M(\r)$ with $2r+2$ corners, and is
homotopic to $\a$.  Note that from the construction all the corners
are in the hyperbolic part of $M(\r)$.

By Proposition 2.2, $\b$ bounds a surface $F$ of disk type in $M(\r)$,
such that $$\area(F \cap N^b) < \area(F \cap (M_0 \cup N^b)) \leq
(2r+2)\pi - 2\pi = 2r\pi. \tag{4-2}$$ After a small perturbation rel
$\bdd$ we may assume that $F$ is transverse to $T$.  Let $A = F \cap
T$.  Since $\bdd A = \bdd F \cap T = \cup \bdd \b_i$, $A$ has exactly
$p$ arc components.  As in the proof of Theorem 4.1, the minimality of
$p$ implies that each arc $a_i$ of $A$ is outmost on $F$ in the sense
that it cuts off a disk $\Delta_i$ with interior containing no arc
components of $A$.  We can label $a_i$ and $\Delta_i$ such that either
$\bdd \Delta_i = a_i \cup \b_{2i-1}$ for all $i$, or $\bdd \Delta_i =
a_i \cup \b_{2i}$ for all $i$.

If $\Delta_i \cap \bdd F = \b_{2i-1}$ for all $i$, then since the
geodesic arc $\b_{2i-1}$ in $M_0$ cannot be homotopic in $M$ to the arc
$a_i$ on $T$, there must be some circle component $\mu_i$ of $A$ in $\Int
\Delta_i$ which is nontrivial on $T$. Applying Lemma 4.3  to the disks
$B_i$ in $\Delta_i$ bounded by $\mu_i$, we get
$$ \area(F\cap N^b) \geq \sum \area(B_i \cap N^b) \geq p(12\pi -
(2\pi+1)) > 2r\pi$$
which is a contradiction to (4-2).

Now assume $\bdd \Delta_i = a_i \cup \b_{2i}$ for all $i$.  Consider a
$\Delta_i$ such that $\b_{2i} \neq \a_{2i}$, $i<p$.  Recall from the
definition of $\b_i$ that there are exactly $r$ such arcs.  We have
shown that in this case $t(\a_{2i}) > 3\pi$ and $b_i > 5$, and we want
to show that $\area(\Delta_i\cap N^b) \geq 2\pi$.  This follows from
Lemma 4.2 if some circle component of $A$ in $\Delta_i$ is nontrivial
on $T$.  Hence assume that all circle components of $A$ in $\Delta$
are trivial on $T$.  In particular, if we denote by $Q$ the component
of $\Delta_i$ cut along $A$ which contains $\bdd \Delta_i$, then $\bdd
Q - \bdd \Delta_i$ is null-homotopic on $T_i$, so $\bdd \Delta_i =
\a_i\cup \b_{2i}$ is also null-homotopic in $N_i(\r_i)$.  Let $\b'$ be
an arc on $T_i$ homotopic to $\b_{2i}$ in $N_i$.  Then $\b'\cup a_i$
is null-homotopic in $N_i(\r_i)$, hence $[(\bdd Q \cap T_i)\cup \b'] =
[a_i\cup \b'] = k[\r_i] \in H_1(T_i)$ for some $k$.  We can now apply
Lemma 3.4 to the surface $Q = \Delta_i \cap N_i^b$: If $k = 0$ then
$$\area(Q) > \area(R_2(t(\b'), b_i)) > 2\pi.$$
The last inequality follows from Lemma 3.1(3) because
we have shown that $t(\b') = t(\a_{2i})>3\pi$
and $b_i > 5$.  If $k \neq 0$, choose a slope $\delta$ as in the statement of the theorem.  Then by (4-1) and Lemma 4.3 we have
$$\align \area(Q) & \geq t_{\d}(k\r_i)(1 - \frac 1b) - t_{\d}(\a) \geq
t_{\d}(\r_i) - t_{\d}(\a) - \frac {t_{\d}(\r_i)}b \\ & \geq
t_{\d}(\r_i) - t_{\d}(\a) - \frac {t(\r_i)}b \geq 5\pi - (2\pi + 1) >
2\pi.  \endalign$$ In either case, $\area(\Delta_i \cap N^b) \geq
\area(Q) > 2\pi$.  Since there are exactly $r$ outmost disks
$\Delta_i$ with $\b_{2i} \neq \a_{2i}$, it follows that $\area(F\cap
N^b) \geq 2r\pi$, which is again a contradiction to (4-2).  \endproof

\head 5. Dehn surgery and Freedman tubing of immersed surfaces
\endhead

Let $M$ be a complete hyperbolic 3-manifold.  For $\mu$ a small
positive number, let $N = N_{\mu}$ be the toroidal cusp components of
$M_{(0,\mu]}$, and $T = T_{\mu} = \bdd N_{\mu}$.  Let $M_0 = M - \Int
N$.  Then $M = N \cup _T M_0$.

A $\pi_1$-injective surface $F: S \to M$ is {\it geometrically
finite\/} if $F_* (\pi_1 S_i)$ is a geometrically finite subgroup of
$\pi_1 M \subset PSL_2(\Bbb C)$ for each component $S_i$ of $S$.  We
need some basic facts about geometrically finite surface groups.  One
is referred to [Mg] for more details.

Assume that $S$ is connected, and $F: S \to M$ is a hyperbolic,
geometrically finite surface in a complete hyperbolic manifold $M$.
Consider the covering $p: X = X_{\Gamma} \to M$ corresponding to the
subgroup $\Gamma = F_*(\pi_1 S)$ of $\pi_1(M)$.  Then $X$ is a
geometrically finite complete hyperbolic manifold.  Denote by $C(F) =
C(X)$ the convex core of $X$, which by definition is the quotient
$C_{\Gamma}/\Gamma$, where $C_{\Gamma}$ is the convex hull of the
limit set of $\Gamma$, and the action of $\Gamma$ on $C_{\Gamma}$ is
induced by the action of $\Gamma$ on its limit set.  Since $\Gamma$
contain no $\Bbb Z^2$ subgroup, the following is a special case of
Lemma 6.5 and Theorem 6.6 of [Mg].

\proclaim{Lemma 5.1} There is an $\e_0 > 0$ such that if $0< \e \leq
\e_0$, then (i) $C(X) \cap X_{[\e, \infty)}$ is compact, (ii) $X_{(0,
\e]}$ has only finitely many components, and (iii) each component of
$X_{(0,\e]}$ is a $\Bbb Z$-cusp, which intersects $C(X)$ in a set
isometric to
$$ \{ (x,y,z) \in \H^3 \,| \, z \geq 1 \text{ and } A_1 \leq y \leq
A_2 \} / (g),$$ where $g$ is a translation in the $x$-direction, and
$A_1, A_2$ are constant depending on the cusp.
\qed \endproclaim

The lifting of $N = N_{\mu}$ to
$X$ is a set of horoballs and $\Bbb Z$-cusps.  Denote by $\tilde N$
the $\Z$-cusp components of $p^{-1}(N)$, and let $\tilde T = \bdd
\tilde N$.  When $\mu$ is small enough, each component of $\tilde
N$ is a component of $X_{(0,\e]}$ for some $\e \leq \e_0$, so we can
define $\mu(F)$ to be the maximum $\mu$ such that this property
holds.  Below we will always assume that $N = N_{\mu}$ and $T =
T_{\mu}$ has been chosen such that $\mu = \mu(F)$.  Note that
we usually assume that $F$ has boundary on $T$.  When we rechoose $T =
T_{\mu(F)}$, we extend $F$ add some collars to $\bdd F$ so that
$\bdd F$ still lies in $T$.  Since $\mu(F)$ depends only on the
group $\Gamma = F_*(\pi_1 S)$, this will not cause a logic problem.

Let $P = \tilde T \cap C(X)$.  By Lemma 5.1, $P$ is a finite set of
compact annuli, one for each component of $\tilde T$.  The {\it
width\/} of a component $P_i$ of $P$ is defined as $w(P_i) = A_2 -
A_1$, where $A_i$ are as in Lemma 5.1.  Define $w(F)$ to be the
maximum of $w(P_i)$ over all component $P_i$ of $P$.  (If $F$ is
disconnected, take the maximum over all $P$ corresponding to all
components of $F$.)

The core of $P_i$ projects to a curve $\a_i'$ on $T$, which is a
nontrivial multiple of some slope $\a_i$ on $T$, usually called a
parabolic slope of $F$.  Since $\pi_1 X = \pi_1 F$, $\a'_i$ is
homotopic to a nontrivial curve on $F$, hence a parabolic slope is a
coannular slope.  The reverse is also true: If a nontrivial curve
$\a'_i$ on $T$ is homotopic to a curve on $F$, then it represents a
parabolic element of $\pi_1 F$, so its lifting on $X$ is homotopic
into some $\Bbb Z$-cusp, hence homotopic to some nontrivial curve on
some $P_i$.  Therefore, the set of parabolic slope of $F$ are the same
as the set of coannular slopes of $F$ on $T$.  By Lemma 5.1, $T$ has
only finitely many coannular slopes of $F$.  The following theorem
says that if the Dehn filling slope is far away from all coannular
slopes of $F$, then $F$ remains $\pi_1$-injective after Dehn filling.

\proclaim{Theorem 5.2} Let $F$ be a hyperbolic, geometrically finite
surface in $M$.  Let $\r = (\r_1, ..., \r_n)$ be a multiple slope on
$T$ such that $t(\r_i) \geq 12\pi$ and $t_{\b}(\r_i) \geq w(F) + 5\pi$
for all coannular slopes $\b$ of $F$.  Then $F$ is $\pi_1$-injective
in $M(\r)$.  \endproclaim

\proof We need to show that if $\a$ is a nontrivial curve on $F$, then
it is also a nontrivial curve in $M(\r)$.  Let $\tilde \a$ be its
lifting to $X= X_F$.  Then $\tilde \a$ is homotopic to a geodesic
$\tilde \a'$ in the convex hull $C(X)$.  The intersection of $\tilde
\a$ with $\tilde T$ cuts $\tilde \a$ into arcs $\tilde \a_1, ...,
\tilde \a_{2n}$, where $\tilde \a_{2i-1}$ lies in $X_{[\e, \infty)}$,
and $\tilde \a_{2i}$ on the cusps.  By the choice of $T =
T_{\mu(F)}$, the image of $C(X) \cap X_{[\e, \infty)}$ is disjoint
from the interior of $N$, hence the projection of $\tilde \a_i$ gives
a decomposition $\a = \a_1 \cup ... \cup \a_{2n}$, with $\a_{2i}$ the
components of $\a\cap N$.  Each $\tilde \a_{2i}$ is homotopic to an
arc lying on a strip of width at most $w(F)$ bounded by geodesics
homotopic to the lifting of some coannular slope $\b$ of $T$, hence
$t_{\b}(\a_{2i}) \leq w(F)$.  By assumption we have $t_{\b}(\r_i) -
t_{\b}(\a_{2i}) \geq t_{\b}(\r_i) - w(F) \geq 5\pi$.  Therefore by
Theorem 4.3 $\a$ is a nontrivial curve in $M(\r)$.  \endproof

The most interesting case is when $F$ is a closed essential surface in
a compact hyperbolic manifold $W$.  The following theorem says that
when finitely many strips centered at coannular slopes and finitely
many other slopes are excluded from the space of Dehn filling slopes,
then $F$ survives surgery.  Note that $F$ is not assumed to be
geometrically finite.

\proclaim{Theorem 5.3} Let $T$ be a set of tori on the boundary of a
compact, orientable, hyperbolic 3-manifold $W$.  Let $F$ be a compact
essential surface in $W$ with $\bdd F \subset \bdd M - T$, and let
$\b$ be the set of coannular slopes of $F$ on $T$.  Then there is an
integer $K$ and a finite set of slopes $\Lambda$ on $T$, such that $F$
is $\pi_1$-injective in $W(\r)$ for all multiple slopes $\r$ on $T$
satisfying $\Delta(\r, \b) \geq K$ and $\r_i \notin \Lambda$.
\endproclaim

\proof We first assume that $\bdd W$ is a set of tori.  Since $W$ is
hyperbolic and $F$ is essential, no component of $F$ is an annulus or
torus, hence $F$ is hyperbolic.  Let $M$ be the interior of $W$, which
by definition has a complete hyperbolic structure.  Since $F$ is
disjoint from $T$, it cannot be a virtual fiber, hence according to
Bonahon and Thurston [B,Th1] it is automatically geometrically finite.
More explicitly, assume that $F$ is geometrically infinite and let
$X_F$ be the covering of $M$ corresponding to the subgroup $\pi_1(F)$.
Then Bonahon [B] showed that every end of $X_F$ relative to the cusp
neighborhoods is geometrically tame, while Thurston [Th1, Theorem
9.2.1] showed that every end of $X_F$ relative to cusp neighborhoods
which is geometrically tame and geometrically infinite must either
correspond to a virtual fiber for $M$ or project to a geometrically
tame and geometrically infinite end of $M$ modulo cusp neighborhoods.
Since we have assumed that $\bdd W$ is a set of tori, $M$ has no
geometrically infinite end modulo cusp neighborhoods, therefore $F$
must be a virtual fiber, which is a contradiction.

Identify the manifold $M_0$ above with $W$, so $\bdd W = T =
T_{\mu(F)}$.  Let $\Lambda$ be the set of slopes $\lambda$ on $T$ such
that $t(\lambda) < 12\pi$.  For each $\b_i$ on some $T_j$, define
$u(\b_i) = \area(T_j) / t(\b_i)$.  Then $t_{\b_i}(\r_j) = \Delta(\b_i,
\r_j) u(\b_i)$.  Choose $K$ so that $K > (w(F) + 5\pi)/u(\b_i)$ for
all $i$.  The result then follows from Theorem 5.2.

Now assume that $W$ has some higher genus boundary components.  If
$\bdd W$ is compressible, then by an innermost circle outermost arc
argument one can show that $F$ can be homotoped to be disjoint from a
maximal set of compressing disks $D$.  Let $W'$ be $W$ cut along $D$.
($W' = W$ if $D = \emptyset$.)  Then $F$ is essential in $W'$ except
that it is possibly homotopic to a surface in a non-torus component of
$\bdd W'$.  Let $\hat W'$ be the double of $W'$ along the non-torus
components of $\bdd W'$.  Denote by $\hat F,\hat T, \hat \b,
\hat \r$ the
double of $F, T, \b, \r$ in $\hat W'$, respectively.  By an innermost
circle outermost arc argument one can show that $\hat F$ is
$\pi_1$-injective in $\hat W$.  Let $q: \hat W' \to W'$ be the obvious
quotient map.  If $A$ is an annulus in $\hat W'$ with one boundary
component on each of $\hat F$ and $\hat T$, then $q(A)$ is an annulus
in $W'$ with one boundary component on each of $F$ and $T$.  Hence
$\hat \b$ is the set of all coannular slopes of $\hat F$ in $\hat W'$.
By the above, there is a number $K$ and a set of slopes $\Lambda'$
such that $\hat F$ is $\pi_1$-injective in $\hat W'(\hat \r)$ when
$\Delta(\hat \r, \hat \b)
\geq K$ and $\hat \r_i \notin \Lambda'$.  Since $F$ is $\pi_1$-injective in
$\hat F$, we see that $F$ is $\pi_1$-injective in
$W'(\r)$.  Since $W(\r)$ is obtained from $W'(\r)$ by adding some
1-handles, $F$ is also $\pi_1$-injective in $W(\r)$.  Let $\Lambda
= q(\Lambda')$.  Then the result follows.
\endproof

We now consider Freedman tubings of essential surfaces.  Let $\hat S$
be a surface containing $S$, such that $\hat S - S$ is a set
of annuli.  Then a surface $\hat F: \hat S \to M_0$ is called a {\it
Freedman tubing\/} of $F$ if $\hat F|_S = F$, and $\hat F(\hat S - S)
\subset T$.  We will use $A = \hat F - \Int F$ to denote the
restriction of $\hat F$ to $\hat S - \Int S$, and call a component
$A_i$ of $A$ a {\it tubing annulus.}  Let $\d_i$ be a component of
$\bdd A_i$.  Then the {\it length\/} of a tube $A_i$ is defined as
$$\ell(A_i) = \min \, \{ \, t_{\d_i}(\a) \, | \, \a \text{ an
essential arc on $A_i$}\}$$ Denote by $\ell(A) = \min \ell(A_i)$.
Clearly, $\ell(A_i)$ would become very large when $A_i$ wraps around
the torus many times.  For example, if $A_i\subset T_j$ is immersed
and contains a sub-annulus $A'_i$ with both boundary components on the
same geodesic curve of $T_j$, and $A'_i$ wraps $k$ times around $T_j$,
then $\ell(A_i) \geq k \, \area(T_j) / t(\d_i)$.

Theorem 5.7 below says that a Freedman tubing of a geometrically
finite surface is essential if the tubes are long enough.  This
generalizes a result of Freedman-Freedman [FF] and Cooper-Long [CL2]
(see also [Li]), where the above result is proved for embedded
surfaces.  In most cases, one can apply Theorem 5.2 to show that it
remains essential after most Dehn fillings.  The assumption that $F$
is geometrically finite is necessary: if $F$ is geometrically
infinite, then $F$ is a virtual fiber, hence all Freedman tubings of
$F$ are inessential.

A boundary component $\d_i$ of $F$ can be pushed around $T$ many
times.  We need a number to measure how far $\d_i$ is away from a
standard position.  We would consider $F$ to be in a standard position
if its lifting to $X$ lies in the convex core $C(X)$.  Let $\tilde
\d_i$ be the component of $\bdd \tilde F$ which projects to $\d_i$.
Each $\tilde \d_i$ is on some component $\tilde T_i$ of $\tilde T$,
which contains a component $P_i$ of $P$.  Define a number $\rho(\d_i)$
to be the minimum nonnegative number such that $\tilde \d_i$ lies in a
$\rho(\d_i)$ neighborhood of $P_i$ on $\tilde T_i$.  Since $\tilde
\d_i$ is compact on $\tilde T_i$, such a number exists.  If $F$ is a
(possibly disconnected) geometrically finite surface in $M$ with some
boundary components on $T = T_{\mu(F)}$, define
$$\rho(F) = \max \, \rho(\d_i)$$ where the maximum is taken over all
boundary components of $F$ which is to be tubed.

\proclaim{Lemma 5.4} Let $\a$ be an arc on $\tilde T$ with one
endpoint $p_1$ on $\tilde \d_i$ and the other endpoint on $P_i$.
Then
$$ t_{\tilde \d_i}(\a) \leq w(F) + \rho(F).$$
\endproclaim

\proof Homotope $\a$ to $\a_1 \cdot \a_2$, where $\a_1$ is a shortest
arc from $p_1$ to some point in $P_i$, and $\a_2$ an arc in $P_i$.
Since $P_i$ is a strip bounded by geodesics of $\tilde T_i$ parallel
to $\tilde \d_i$, by definition we have $t_{\tilde \d_i}(\a_1)
\leq w(F)$, and $t_{\tilde \d_i}(\a_2) \leq \rho(F)$.  \endproof

Two arcs $\a_1, \a_2$ in $X$ with $\bdd \a_i \subset \tilde T$ are
{\it $\tilde T$-homotopic\/} if there are arcs $\b', \b''$ on $\tilde
T$ such that $\a_1 \sim \b'\cdot \a_2 \cdot \b''$.  Clearly this is an
equivalence relation.  An arc $\a$ in $X$ is of type I if it projects
to a type I arc in $M$.

\proclaim{Lemma 5.5} Any proper essential arc $\a$ of $\tilde F$ is
$\tilde T$-homotopic to a type I arc of $X$ with endpoints on $P$.
\endproclaim

\proof First deform $\a$ by a $\tilde T$-homotopy to an arc $\a_1$
with $\bdd \a_1 \subset P$.  This is possible because each component
of $\tilde T$ contains a component of $P$.  Now homotope $\a_1$ (rel
$\bdd$) to a geodesic $\a_2$ in $X$.  Since $C(X)$ is a convex set,
$\a_2 \subset C(X)$, so $\a_2 = \b_1 \cdot \a_3 \cdot \b_2$, where
$\a_3$ is a geodesic of type I with endpoints in $P$, and $\b_1, \b_2$
are (possibly empty) arcs in $C(X) \cap X_{(0, \e]}$, which can be
pushed into $\tilde T$, hence $\a_2$ is $\tilde T$-homotopic to
$\a_3$.  \endproof

\proclaim{Lemma 5.6} Let $F$ be a geometrically finite surface in $M$.
Let $\a$ be an essential arc of $F$ with endpoints on boundary
components $\d_0, \d_1$ of $F$ which lie on $T = T_{\mu(F)}$.  Then
$\a$ is homotopic to $\b_0 \cdot \a' \cdot \b_1$, where $\a'$ is an
arc of type I, and $\b_i$ are arcs on $T$ with $t_{\d_i}(\b_i) \leq
\rho(F) + w(F)$.  \endproclaim

\proof Consider the lifting $\tilde a$ of $\a$ on $\tilde F \subset
X$.  By Lemma 5.5, $\tilde \a$ is homotopic to $\tilde \b_0 \cdot
\tilde \a' \cdot \tilde \b_1$, where $\tilde \a'$ is of type I, and
$\tilde \b_i$ is an arc on some component $\tilde T_i$ of $\tilde T$
with one endpoint on each of $\tilde \e_i$ and $P_i$.  Projecting
these curves into $M$, we get $\a \sim \b_0 \cdot \a' \cdot \b_1$.  By
Lemma 5.4, we have $ t_{\d_i}(\b_i) = t_{\tilde \d_i} (\tilde
\b_i) \leq \rho(F) + w(F).$ \endproof

Recall that the wrapping number of an annulus $A$ on a torus $T$ is
defined as
$$ \text{wrap}(A) = \{\, |A\cdot p| \, \, | \, p \in T\, \}$$
where $A\cdot p$ denotes the algebraic intersection number between $A$
and $p$, which is well defined for all points $p\notin \bdd A$.

\proclaim{Theorem 5.7} Let $F$ be a geometrically finite surface in a
compact hyperbolic 3-manifold $W$.  Then there is a constant $K$ such
that if $\hat F$ is a Freedman tubing of $F$ with $\text{wrap}(\hat F,
F) \geq K$, then $\hat F$ is $\pi_1$-injective in $W$.
\endproclaim

\proof
Let $M$ be the interior of $W$.  By assumption $M$ is a complete
hyperbolic manifold.  Identify $M_0$ above with $W$, possibly with
some higher genus boundary components removed.  Let $T = \bdd M_0$.
Clearly $\ell(\hat F - \Int F)$ goes to infinity when
$\text{wrap}(\hat F, F)$ approaches infinity.  Choose $K$ large enough
such that when $\text{wrap}(\hat F, F) > K$, we have $\ell(\hat F -
\Int F) > 2 (\rho(F) + w(F) + \pi)$.

We need to show that any nontrivial curve $\a$ on $\hat F$ is also
nontrivial in $M$.  If $\a$ is homotopic to a curve on $F$ or $A =
\hat F - \Int F$ then $\a$ is nontrivial in $M$ because $F$ is $\pi_1$
injective.  So assume $\a = \a_1 \cup ... \cup \a_{2n}$, where
$\a_{2i-1} \subset F$ and $\a_{2i} \subset A$ are essential arcs.  By
Lemma 5.6, we have $\a_{2i-1} \sim \b_{2i-1} \cdot \a'_{2i-1} \cdot
\r_{2i-1}$, where $\a'_{2i-1}$ is a type I arc, and
$t_{\d'}(\b_{2i-1})$ and $t_{\d''}(\r_{2i-1}) \leq \rho(F)+w(F)$,
where $\d', \d''$ are boundary components of $F$ containing the
endpoints of $\a_{2i-1}$.  Put $\a'_{2i} = \r_{2i-1} \cdot \a_{2i}
\cdot \b_{2i+1}$.  Then $\a
\sim \a'_1 \cdot \a'_2 \cdot \cdots \cdot \a'_{2n}$. Let $\d_i$ be the
boundary component of $F$ containing an endpoint of $\a_{2i}$.  Then
$$ \align t(\a_{2i}') & \geq t_{\d_i}(\a'_{2i}) \geq
t_{\d_i}(\a_{2i}) - t_{\d_i}(\r_{2i-1}) - t_{\d_i}(\b_{2i+1}) \\
& \geq \ell(\hat F - \Int F) - 2 \rho(F)-2w(F) > 2\pi \endalign$$
Therefore by Theorem 4.1, $\a$ is a nontrivial curve in $M$.
\endproof

\head 6.  Upper bounds on surgery distance and tubing length.
\endhead

Theorems 5.3 is the best possible in the sense that there is no
universal bounds on the number $K$ in the theorem.  Similarly, Theorem
5.7 is the best possible in the sense that there is no universal bound
on how many time a surface need to tube around a torus boundary
component in order to produce an essential surface.  Assume that $\hat
F$ is a Freedman tubing of an essential surface $F$, with tubes on a
torus $T = \bdd M_0$.

\proclaim{Theorem 6.1} (i) For any constant $K$, there is an embedded,
geometrically finite surface $F$ in a hyperbolic manifold $M$, such
that all Freedman tubing $\hat F$ of $F$ with $\text{wrap}(\hat F)
\leq K$ are inessential.

(ii) For any constant $K$, there is a closed essential surface $F$ in
a hyperbolic manifold $M$, and a slope $\b$ on $T$, such that $F'$ is
compressible in $M(\r)$ for all $\r$ with $\Delta(\r, \b) \leq K$.
\endproclaim

\proof (1) Let $S$ be a compact orientable surface of genus $g > K$
with a single boundary component $c$.  Let $\a_1, ..., \a_g$ be a set
of mutually disjoint nonseparating curves cutting $S$ into a connected
planar surface.  By Theorem 1.1 of [WWZ], there exists a pseudo-Anosov
map $\varphi: S \to S$ such that $\varphi(\a_i) = \a_{i+1}$ for $i<g$.
(Note that $\varphi(\a_g) \neq \a_1$, otherwise $\varphi$ would be
reducible.)  Let $W = S \times I - N(\a'_1)$, where $\a'_1$ is the
curve $\a_1 \times \frac 12$ in the interior of $S\times I$ isotopic
to $\a_1$.  Let $M = W/((x,1) \sim (\varphi(x),0))$.  Since $\varphi$
is pseudo-Anosov, it is easy to check that $M$ is irreducible and
atoroidal, and it cannot be a Seifert fiber space because $S = S\times
0$ is an essential hyperbolic surface in $M$ disjoint from one
boundary component of $M$.  Therefore by Thurston's hyperbolization
theorem for Haken manifolds [Th2], $M$ is hyperbolic.

Let $F$ be the disjoint union of two copies of $S$ with opposite
orientation.  Then $F$ is $\pi_1$ injective, and is not a virtual
fiber because it is disjoint from one boundary component of $M$.
Hence it is geometrically finite.  Let $\hat F$ be a Freedman tubing
of $F$ with $\text{wrap}(\hat F) = k \leq K$.  We want to show that
$\hat F$ is inessential in $M$.

Let $\tilde M$ be the infinite cyclic covering of $M$ dual to the
surface $S$.  Note that $\tilde M$ can be constructed by taking
infinitely many copies of $W$, denoted by $W_i$ ($i\in \Bbb Z$), and
gluing the surface $S \times 1$ in $W_i$ to $S\times 0$ in $W_{i+1}$
using the map $\varphi$.  Let $X_k$ be the union of $W_1, ..., W_k$ in
$\tilde M$.  Then $\hat F$ lifts to a surface in $\tilde M$ homotopic
to $\bdd X_k$.  Put $\a'_i = \a_i \times \frac 12$.  One can check
that when $k\leq K$, $X_k$ is homeomorphic to the manifold $(S \times I)
- \a'_1\cup ... \cup \a'_k$.  Let $\b$ be an essential arc on $S$
disjoint from all $\a_i$.  Then $\b \times I$ is a compressing disk of
$\bdd X_k$.  It follows that $\hat F$ is compressible in $M$.

(2) Let $\hat F$ be a Freedman tubing of $F$ such that $\hat F$ is
essential in $M$, and the wrapping number $w$ of $\hat F$ is minimal
among all such surface.  Since $F$ is geometrically finite and
embedded, the existence of such a surface follows from [CL2] or [Li],
or from Theorem 5.7.  Let $\b$ be the boundary slope of $F$, and let
$A$ be the tubing annulus $\hat F - \Int F$.  Assume that $\Delta =
\Delta(\r, \b)\leq K$.  Notice that the annulus $A$ is rel $\bdd$
homotopic in the Dehn filling solid torus to another annulus $A'$ on
$\bdd M$ with wrapping number $w' = |w - \Delta|$, so $\hat F$ is
homotopic in $M(\r)$ to a surface $\hat F' = F \cup A'$ which is a
Freedman tubing of $F$ with wrapping number $w'$.  By the choice of
$w$, it follows that $\hat F$ is inessential in $M(\r)$ for all $\r$
such that $\Delta(\r, \b) < 2w$.  By (1) we have $2w > w >K$, hence
the result follows.  \endproof

Although there is no universal upper bound on the wrapping number of
an essential Freedman tubing surface, it has been shown by Li [Li] that
an upper bound in terms of genus and number of boundary components of
$F$ does exist if $F$ is an embedded surface.  Li showed that if $F$
is embedded with genus $g$ and $b$ boundary components, then a
Freedman tubing of $F$ is essential if its wrapping number is at least
$6g+2b-3$.

\proclaim{Problem 6.2} Find the minimal constant $C(g,b)$ such that
if $F$ is a geometrically finite {\it embedded} surface with genus $g$
and $b$ boundary components, then all Freedman tubing of $F$ with
wrapping number at least $C(g,b)$ is essential.
\endproclaim

Li's result [Li] shows that $C(g,b) \leq 6g+2b-3$, and the proof of
Theorem 6.1 shows that $C(g,b) > g$.

For immersed surface, no such number would exist if we do not assume
that $F$ is in standard position.  The reason is because we can slide
one component of $\bdd F$ around the torus many times, so when tubing
on the opposite direction, a long part of the tube would just homotope
that boundary component of $F$ back to its original position.
However, one can consider the number of tubings which is inessential.
For the embedded case, it is at most $2 C(g,b) + 1$.  For simplicity
let us consider the case that $F$ has only two boundary components.

\proclaim{Conjecture 6.3} Let $F$ be a surface with two boundary
components, both on a torus component of $\bdd M$.  Let $g$ be the
genus and $b$ the number of boundary components of $F$.  Then exists a
constant $C'(g,b)$ depending only on $g$ and $b$, such that all but at
most $C'(g, b)$ of the Freedman tubings of $F$ are $\pi_1$-injective.
\endproclaim

The following result gives an estimation of tubing length when $F$ is
a totally geodesic surface, which leads to an upper bound on wrapping
numbers of inessential Freedman tubing in this special case.
Existence of immersed totally geodesic surfaces can be found in [AR]
and [Re].

\proclaim{Theorem 6.4} Let $T$ be the boundary tori of a set of
disjoint cusps $N$ in $M$, let $F'$ be a totally geodesic surface in
$M$, and let $F = F' \cap M_0 = M - \Int N$.  If $\hat F$ is a
Freedman tubing of $F$ with $\ell(\hat F - \Int F) \geq \pi$, then
$\hat F$ is $\pi_1$-injective in $M$.  \endproclaim

\proof Notice that we do not require $T = T_{\mu(F)}$.  The
intersection of $F'$ with $N$ is a set of totally geodesic annuli,
hence they are perpendicular to $T = \bdd N$.  A nontrivial curve $\a$
on $\hat F$ can be homotoped on $\hat F$ either to a curve on $F$ or
to a curve $\a_1 \cup ... \cup \a_{2n}$ with $\a_{2i-1}$ a geodesic
arc on $F$ perpendicular to $T$, and $\a_{2i}$ an essential arc on
$\hat F - F$.  Since $F$ is totally geodesic, $\a_{2i-1}$ is also a
geodesic of $M$, hence the result follows from Theorem 4.1.  \endproof

\proclaim{Corollary 6.5} Let $F$ be as in Theorem 6.4.  If $\hat F$ is
an inessential Freedman tubing of $F$, then $\text{wrap}(\hat F) \leq
2\pi^2(2g+b-2)/\sqrt 3$.  \endproclaim

\proof Let $A$ be a tubing annulus of $\hat F$, and let $\b$ be the
boundary slope of $A$.  Extend $F$ to a complete hyperbolic surface
$F'$ by adding a cusp at each of its boundary component.  By the
Gauss-Bonnet theorem, $\area(F') = 2\pi(2g+b-2)$.  Each cusp with
boundary on $\bdd A$ has area $= t(\b)$, and there are two of them,
hence $t(\b) < \pi(2g+b-2)$.  Choose $T$ to be a set of maximal cusps,
then each slope of $T$ has length at least $1$, hence the area of each
component of $T$ is at least $\sqrt 3/2$.  If $\r$ and $\b$ are on the
torus $T_i$, then
$$\ell(A) \geq \text{wrap}(A) \area(T_i)/t(\b) \geq \text{wrap}(\hat
F) \sqrt 3 / 2\pi(2g+b -2).$$ By Theorem 6.4, $\hat F$ is essential if
$\ell(A) \geq \pi$ for all tubes $A$ of $\hat F$, which is true if
$\text{wrap}(\hat F) \geq 2\pi^2(2g+b-2)/\sqrt 3$.  \endproof

If $F$ is a closed, embedded, incompressible surface in $M_0$ which is
not coannular to a torus $T \subset \bdd M_0$, then Theorem 1 of [Wu]
says that $F$ remains incompressible in $M(\r)$ for all but at most
three slopes $\r$ on $T$.  For immersed essential surfaces $F$ in $M$
without coannular slopes (also called accidental parabolics), Theorem
1.1 says that $F$ remains essential in $M(\r)$ except for finitely
many $\r$ on $T$.  The answer to the following problem is likely to be
negative.  If $M$ is not assumed hyperbolic, there are examples
showing that no upper bound exists.  However, no examples are known
for hyperbolic manifolds.

\proclaim{Problem 6.6} Let $F$ be a closed essential surface in a
hyperbolic manifold $M$, and assume that $F$ has no coannular slopes.
Does there exist a universal upper bound on the number of slopes $\r$
on a torus boundary component $T$ of $M_0$ such that $F$ is
inessential in $M(\r)$?  \endproclaim

Many hyperbolic manifolds do not contain closed {\it embedded\/}
essential surfaces.  However, it was proved in [CLR] that any
hyperbolic $M$ with some toroidal cusps contains a closed essential
surface.  The surfaces constructed there are Freedman tubings of some
surfaces in certain covering spaces of $M$, and hence all have
coannular slopes.  The following seems to be an interesting open
problem.  The corresponding problem for closed hyperbolic manifold is
also open, and is part of the virtual Haken conjecture.

\proclaim{Conjecture 6.7} Every hyperbolic manifold with toroidal
cusps contains a closed essential surface without coannular slopes.
\endproclaim

\enddocument